\documentclass[11pt]{article}
\usepackage{mathrsfs}
\usepackage{epic,eepic,epsf,epsfig}
\usepackage{amsfonts,srcltx,mathrsfs}
\usepackage{amsmath}
\usepackage{amssymb}
\usepackage{amsbsy}
\usepackage{graphicx}
\usepackage{amsfonts}
\usepackage{color}
\usepackage{mathrsfs}
\usepackage{epic,eepic,epsf,epsfig}
\usepackage{graphicx}
\usepackage{amsfonts,srcltx,mathrsfs}
\usepackage{array}
\usepackage{amsthm}
\newtheorem{theorem}{Theorem}
\newtheorem{lem}{Lemma}

\newtheorem{clm}{Claim}

\def\f{\noindent}

\def\demo{\f {\bf Proof.}\hskip10pt}

\begin{document}

\markboth{ et. al}{Cacti with $n$-vertices and $t$ cycles\\ having extremal edge-Wiener index}

\title{Bounds on the edge-Wiener index of cacti\\ with $n$ vertices and $t$ cycles}

\author{Siyan Liu, Rong-Xia Hao\footnote{Corresponding author. Email: 14275011@bjtu.edu.cn (Siyan Liu), rxhao@bjtu.edu.cn (Rong-Xia Hao),  he1046436120@126.com (Shengjie  He)}, Shengjie He\\[0.2cm]
{\em\small Department of Mathematics, Beijing Jiaotong University, Beijing 100044, P.R. China}}

\date{}
\maketitle

The edge-Wiener index $W_e(G)$ of a connected graph $G$ is the sum of distances between all pairs
of edges of $G$. A connected graph $G$ is said to be a cactus if each of its blocks is either a
cycle or an edge. Let $\mathcal{G}_{n,t}$ denote the class of all cacti with $n$ vertices and $t$ cycles. In this paper, the upper bound and lower bound on the edge-Wiener index of graphs in $\mathcal{G}_{n,t}$ are identified and the corresponding extremal graphs are characterized.

\medskip

\f {\em Keywords:}  Cactus; Edge-Wiener index; Upper bound; Lower bound.

\section{Introduction}

Throughout this paper, all graphs we considered are finite, undirected, and simple. Let $G$ be a connected graph with vertex set $V(G)$ and edge set $E(G)$. For a vertex $u \in V(G)$, the $degree$ of $u$, denote by $d_{G}(u)$, is the number of vertices
which are adjacent to $u$.
For a vertex $v \in V(G)$, denote by $N_{G}(u)$, the set of the vertices
which are adjacent to $v$.
Call a vertex $u$ a $pendent\ vertex$ of $G$, if $d_{G}(u)=1$ and call an edge $uv$ a $pendent\ edge$ of $G$, if $d_{G}(u)=1$ or $d_{G}(v)=1$.
By $G-v$ and $G-uv$ we denote the graph obtained from $G$ by deleting a vertex $v \in V(G)$, or an edge $uv \in E(G)$, respectively (This notation is naturally extended if more than one vertex or edge are deleted). Similarly, $G+uv$ is obtained from $G$ by adding an edge $uv \notin E(G)$.
For any two vertices $u, v \in V(G)$, let $d_{G}(u, v)$ denote the distance between $u$ and $v$ in $G$. Denote by $P_n$, $S_n$ and $C_n$ a path, star and cycle on $n$ vertices, respectively. We refer to Bondy and Murty~\cite{Bondy} for notation and terminologies used but not defined here.

The Wiener index is one of the oldest and the most thoroughly studied topological indices. The Wiener index of a graph $G$ is defined as
$$W(G)=\sum\limits_{\{ u, v\} \subseteq V(G) }d_{G}(u, v).$$
The edge-Wiener index is defined as the sum of distances between all pairs of edges, namely as
\[
W_e(G) = \sum_{\{f,g\} \subseteq E(G)}d_{G}(f,g).
\]
where $d_G(f,g)$ denotes the distance between $f$ and $g$ in $G$, and also the distance between the corresponding vertices in the line graph of $G$. Note that for any two distinct edges $f = u_1u_2$ and $g = v_1v_2$ in $E(G)$, the distance between $f$ and $g$ equals $d_{G}(f,g) =  \mbox{ min }\{d(u_i, v_j) : i, j \in \{1,2\}\} + 1$. In the case, when $f$ and $g$ coincide, we have $d_G(f,g) = 0$. Nowadays, the Wiener index is a well-known and much studied graph invariant e.g.,~\cite{Do.RI,I.Relation,Martin,Liu,HLP,TAN}.
In ~\cite{EWiener}, Dankelmanna et al. gave bounds on $W_{e}$ in terms of order and size.
But there is not so many conclusions about edge-Wiener index. So we pay attention to it and hope to get more conclusions.

A connected graph $G$ is said to be a $cactus$ if any two of its cycles have at most one common vertex. Let $\mathcal{G}_{n,t}$ be the set of all $n$-vertex cacti, each containing exactly $t$ cycles. A graph is a $chain$ $cactus$ if each block has at most two cut vertices and each cut vertex is shared by exactly two blocks. Obviously, any chain cactus with at least two blocks contains exactly two blocks that have only one cut-vertex. Such blocks are called $terminal$ blocks. If all the cycles in a $cactus$ have exactly one common vertex, then they form a $bundle$. In ~\cite{Ivan}, the Wiener index of cacti with $n$ vertices and $t$ cycles was studied by Gutman. Wang ~\cite{WANGSJ,WANG} determined the the lower bounds on Szeged index and revised Szeged index of cacti with $n$ vertices and $k$ cycles.
He et al.~\cite{S.H.H.Y} determined the lower bounds of edge Szeged index and edge-vertex Szeged index for cacti with order $n$ and $k$ cycles.

In this paper, by using the methods similar to Gutman \cite{Ivan}, the edge-Wiener index of the cacti with $n$ vertices and $t$ cycles is studied.
Moreover, the lower bound on edge-Wiener index of the cacti with given cycles is determined and the corresponding extremal graph is identified. Furthermore, the upper bound on edge-Wiener index of the cacti with given cycles and the corresponding extremal graph are established as well.

\section{Preliminaries}
In this section, we give some preliminary results which will be used in the subsequent sections.

\begin{lem}(~\cite{Martin})\label{lem00}
The Wiener index of $P_n$ and $C_n$ is
$$W(P_n) = \frac{1}{6}n(n + 1)(n - 1),$$
$$ W(C_n)=\left\{
                                                 \begin{array}{ll}
                                                  \frac{1}{8}n(n^2 - 1), & \hbox{if $n$ is odd ;} \\
                                                  \frac{1}{8}n^3, & \hbox{if $n$ is even.}
                                                 \end{array}
                                               \right.
  $$
\end{lem}

\begin{lem}\label{lem-1}
Let $G_1$ and $G_2$ be two connected graphs with disjoint vertex sets where $u_1 \in V(G_1)$, $u_2 \in V(G_2)$. Let $|V(H_i)| = n_i$, $|E(H_i)| = m_i$ for $i = 1,2$, $(n_i \geq 2, m_i \geq 1)$. Construct the graph $G$ by identifying the vertices $u_1$ and $u_2$, and denote the new vertex by $u$ (see Figure 1). Then
\[
W_e(G) = W_e(G_1) + W_e(G_2) + m_1\sum_{f \in E(G_2)}d_{G}(f,u) + m_2\sum_{g \in E(G_1)}d_{G}(g,u) + m_1m_2.
\]
\end{lem}

\begin{figure}
\centering
    \includegraphics[height=4cm,width=6cm]{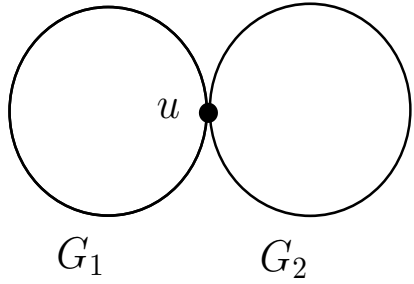}
    \caption{The graph $G$ in Lemma 2}
    \label{fig1}
\end{figure}

\demo  According to the definition of edge-Wiener index, we have
\begin{equation*}
\begin{split}
W_e(G)
      &= \sum_{f,g \in E(G)}d_G(f,g) \\
      &= \sum_{f,g \in E(G_1)}d_{G_1}(f,g) + \sum_{f,g \in E(G_2)}d_{G_2}(f,g) + \sum_{f \in E(G_1), g \in E(G_2)}d_G(f,g) \\
      &= W_e(G_1) + W_e(G_2) + \sum_{f \in E(G_1), g \in E(G_2)}[d_{G}(f,u) + d_G(g,u) + 1] \\
      &= W_e(G_1) + W_e(G_2) + \sum_{f \in E(G_1), g \in E(G_2)}d_G(f,u) + \sum_{f \in E(G_1), g \in E(G_2)}d_G(g,u) + m_1m_2 \\
      &= W_e(G_1) + W_e(G_2) + m_1\sum_{f \in E(G_2)}d_G(f,u) + m_2\sum_{g \in E(G_1)}d_G(g,u) + m_1m_2.
\end{split}
\end{equation*}

That's the end of the proof.
\qed

\begin{lem}\label{lem-2}
Let $G$ be a graph with a cut edge $e'=v_1v_2$, and $G'$ be the graph obtained from $G$ by contracting the edge $e'$ and adding a pendant edge attaching at the contracting vertex; see Figure 2. Let $|E(G_i)|=m_i$ for $i=1,2$. If $d_G(v_i) \geq 2$ for $i=1,2$, we have $W_e(G') < W_e(G)$.

\end{lem}

\begin{figure}
\centering
    \includegraphics[height=4cm,width=10cm]{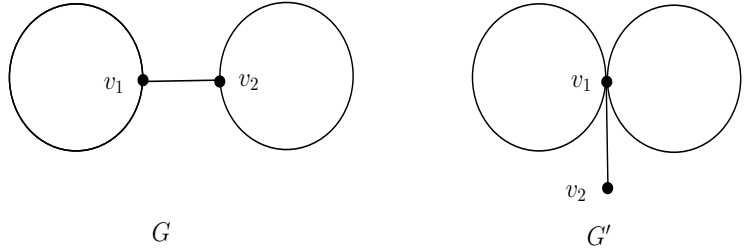}
    \caption{The graphs $G$ and $G'$ in Lemma 3}
    \label{fig2}
\end{figure}

\demo  We calculate the edge-Wiener index of $G$ and $G'$, respectively.
\begin{equation*}
\begin{split}
W_e(G)
      &= \sum_{f,g \in E(G)}d_G(f,g) \\
      &= \sum_{f,g \in E(G_1)}d_{G_1}(f,g) + \sum_{f,g \in E(G_2)}d_{G_2}(f,g) + \sum_{f \in E(G_1), g \in E(G_2)}d_G(f,g) \\
      &\quad + \sum_{f \in E(G_1)}d_G(f, v_1v_2) + \sum_{g \in E(G_2)}d_G(g, v_1v_2) \\
      &= W_e(G_1) + W_e(G_2) + \sum_{f \in E(G_1), g \in E(G_2)}[d_G(f,v_1) + d_G(g,v_2) + 1 + 1] \\
      &\quad + \sum_{f \in E(G_1)}d_G(f, v_1v_2) + \sum_{g \in E(G_2)}d_G(g, v_1v_2).
\end{split}
\end{equation*}
Also, use the same method, then we have
\begin{equation*}
\begin{split}
W_e(G')
      &= \sum_{f,g \in E(G')}d_{G'}(f,g) \\
      &= W_e(G'_1) + W_e(G'_2) + \sum_{f \in E(G'_1), g \in E(G'_2)}[d_{G'}(f,v_1) + d_{G'}(g,v_2) + 1] \\
      &\quad + \sum_{f \in E(G'_1)}d_{G'}(f, v_1v_2) + \sum_{g \in E(G'_2)}d_{G'}(g, v_1v_2).
\end{split}
\end{equation*}
So
\begin{equation*}
\begin{split}
W_e(G) - W_e(G') = \sum_{f \in E(G_1), g \in E(G_2)} 1 = m_1m_2 > 0.
\end{split}
\end{equation*}

That's the end of the proof.
\qed

\begin{lem}\label{lem-3}
Let $G$ be a graph with an even cycle $C_{2k} = v_1v_2 \cdots v_{2k}v_1$ such that $G - E(C_{2k})$ has exactly $2k$ components. Let $G_i$ be the component of $G - E(C_{2k})$ that contains $v_i$ and $|E(G_i)|=m_i$ for $i=1, 2, \cdots, 2k$. Let
\[
G' = G - \sum_{i = 2}^{2k} \sum_{w \in N_{G_i}(v_i)} wv_i + \sum_{i = 2}^{2k} \sum_{w \in N_{G_i}(v_i)} wv_1.
\]
(see Figure 3). Then we have $W_e(G') \leq W_e(G)$ with equality if and only if $C_{2k}$ is an end-block, that is, $G\cong G'$.
\end{lem}

\begin{figure}
\centering
    \includegraphics[height=7cm,width=7cm]{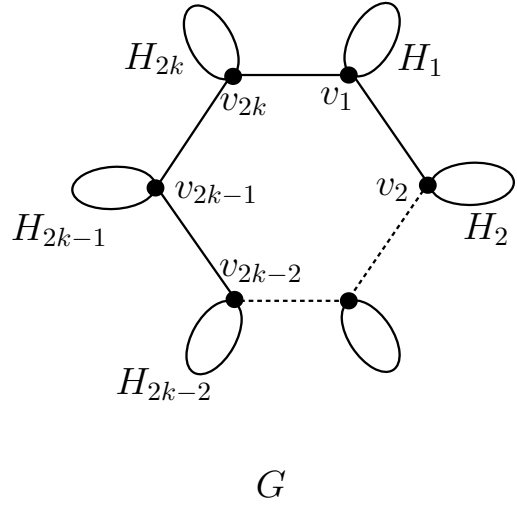}
    \caption{The graph $G$ in Lemma 4}
    \label{fig3}
\end{figure}

\demo
If $C_{2k}$ is an end-block, the lemma holds clearly. Then, one can assume that $C_{2k}$ is not an end-block in  the following.

By the definition of edge-Wiener, we can suppose that $W_e(G) - W_e(G') = T_1 + T_2 + T_3 + T_4$, which are
\begin{equation*}
\begin{split}
T_1 &= \sum_{i = 1}^{2k} \sum_{f,g \in E(H_i)}[d_G(f,g) - d_{G'}(f,g)], \\
T_2 &= \sum_{f,g \in E(C_{2k})}[d_G(f,g) - d_{G'}(f,g)], \\
T_3 &= \sum_{1 \leq i < j \leq 2k} \sum_{f \in E(H_i), g \in E(H_j)}[d_G(f,g) - d_{G'}(f,g)], \\
T_4 &= \sum_{i = 1}^{2k} \sum_{f \in E(H_i), g \in E(C_{2k})}[d_G(f,g) - d_{G'}(f,g)].
\end{split}
\end{equation*}
According to the structure of $G$ and $G'$, we have $T_1 = T_2 = 0$ clearly. Then we calculate $T_3$ and $T_4$.
\begin{equation*}
\begin{split}
T_3 &= \sum_{1 \leq i < j \leq 2k} \sum_{f \in E(H_i), g \in E(H_j)}[d_G(f,g) - d_{G'}(f,g)] \\
    &= \sum_{1 \leq i < j \leq 2k} \sum_{f \in E(H_i), g \in E(H_j)}[d_G(f,v_i) + d_G(g,v_j) + d_G(v_j, v_i) + 1] \\
    &\quad - \sum_{1 \leq i < j \leq 2k} \sum_{f \in E(H_i), g \in E(H_j)}[d_{G'}(f,v_1) + d_{G'}(g,v_1) + 1] \\
    &= \sum_{1 \leq i < j \leq 2k} \sum_{f \in E(H_i), g \in E(H_j)}d_G(v_i, v_j).
\end{split}
\end{equation*}
By $d_G(v_i, v_j) > 0$ when $i \neq j$, we have $T_3 > 0$. Now let's start $T_4$. We say the sum of the edges in each branch $H_i$ is equal to $n_i =\sum_1^{2k} m_i$, each edge of the cycle $C_{2k}$ is $e_i$ for $i = 1,2, \cdots, 2k$. The sum distance between $e_i$ and $g$ ($g \in E(H_i)$) is equal to $\sum_{1 \leq i \leq 2k, 1 \leq j \leq 2k,g \in E(H_j)} d_G(e_i, g)$.

In graph $G$, $d_G(e_i, g) = d_G(g, v_j) + d_G(e_i, v_j) + 1$, and in graph $G'$, $d_{G'}(e_i, g) = d_{G'}(g, v_1) + d_{G'}(e_i, v_1) + 1$, for $v_j$ is the cut vertex of $H_j$. Minus equal parts, we have
\begin{equation*}
\begin{split}
T_4 &= \sum_{i=1}^{2k} \sum_{f \in E(H_i), g \in E(C_{2k})}[d_G(f,g) - d_{G'}(f,g)] \\
    &= \sum_{1 \leq i \leq {2k}, 1 \leq j \leq {2k}, g \in E(H_j)}d_G(e_i, g) - \sum_{1 \leq i \leq {2k}, 1 \leq j \leq {2k}, g \in E(H_j)}d_{G'}(e_i, g). \
\end{split}
\end{equation*}
Now use the previous results,
\begin{equation*}
\begin{split}
T_4 &= \sum_{1 \leq i \leq {2k}, 1 \leq j \leq {2k}}m_jd_G(e_i, v_j) - \sum_{1 \leq i \leq {2k}, 1 \leq j \leq {2k}}m_jd_{G'}(e_i, v_1) \\
    &= \sum_{1 \leq j \leq {2k}}[m_j \sum_{i=1}^{2k} d_G(e_i, v_j)] - \sum_{1 \leq j \leq {2k}}[m_j \sum_{i=1}^{2k}d_{G'}(e_i, v_1)].
\end{split}
\end{equation*}
In graph $G$, the sum distance between the edge $e_i$ and vertex $v_1$ is $\sum_{i=1}^{2k} d_G(e_i, v_1)$ and
\[
\sum_{i=1}^{2k}d_G(e_i, v_1) = [(0 + 1 + 2 + 3 + \cdots k-1) \times 2] = k(k-1).
\]
that means the sum distance of every edge from $C_{2k}$ to $v_1$ is $k(k-1)$. Now let's calculate graph $G'$, the cycle $C_{2k}$ has symmetry, for every cut vertex $v_s$ and $v_t$, we have $\sum_{i=1}^{2k} d_{G'}(e_i, v_s) = \sum_{i=1}^{2k} d_{G'}(e_i, v_t)$, so
\begin{equation*}
\begin{split}
T_4 &= [m_1k(k-1) + m_2k(k-1) + m_3k(k-1) + \cdots + m_{2k}k(k-1)] \\
    &\quad - [m_1k(k-1) + m_2k(k-1) + m_3k(k-1) + \cdots + m_{2k}k(k-1)] \\
    &= 0.
\end{split}
\end{equation*}
Now we know $T_1 = T_2 = T_4 = 0$, so
\[
W_e(G) - W_e(G') = T_1 + T_2 + T_3 + T_4 = T_3 > 0.
\]

That's the end of the proof.
\qed

\begin{lem}\label{lem-4}
Let $G$ be a graph with an odd cycle $C_{2k+1} = v_1v_2 \cdots v_{2k+1}v_1$ such that $G - E(C_{2k+1})$ has exactly $2k+1$ components. Let $G_i$ be the component of $G - E(C_{2k+1})$ that contains $v_i$ and $|E(G_i)|=m_i$ for $i=1, 2, \cdots, 2k+1$. Let
\[
G' = G - \sum_{i = 2}^{2k+1} \sum_{w \in N_{G_i}(v_i)} wv_i + \sum_{i = 2}^{2k+1} \sum_{w \in N_{G_i}(v_i)} wv_1.
\]
Then we have $W_e(G') \leq W_e(G)$ with equality if and only if $C_{2k+1}$ is an end-block, that is, $G\cong G'$.
\end{lem}

\demo
If $C_{2k+1}$ is an end-block, the lemma holds clearly. Then, one can assume that $C_{2k+1}$ is not an end-block in  the following.

By the definition of edge-Wiener, we can suppose that $W_e(G) - W_e(G') = S_1 + S_2 + S_3 + S_4$, where
\begin{equation*}
\begin{split}
S_1 &= \sum_{i = 1}^{2k+1} \sum_{f,g \in E(H_i)}[d_G(f,g) - d_{G'}(f,g)], \\
S_2 &= \sum_{f,g \in E(C_{2k+1})}[d_G(f,g) - d_{G'}(f,g)], \\
S_3 &= \sum_{1 \leq i < j \leq {2k+1}} \sum_{f \in E(H_i), g \in E(H_j)}[d_G(f,g) - d_{G'}(f,g)], \\
S_4 &= \sum_{i = 1}^{2k+1} \sum_{f \in E(H_i), g \in E(C_{2k+1})}[d_G(f,g) - d_{G'}(f,g)].
\end{split}
\end{equation*}
According to the structure of $G$ and $G'$, we have $S_1 = S_2 = 0$ clearly. Then we calculate $S_3$ and $S_4$.
\begin{equation*}
\begin{split}
S_3 &= \sum_{1 \leq i < j \leq {2k+1}} \sum_{f \in E(H_i), g \in E(H_j)}[d_G(f,g) - d_{G'}(f,g)] \\
    &= \sum_{1 \leq i < j \leq {2k+1}} \sum_{f \in E(H_i), g \in E(H_j)}[d_G(f,v_i) + d_G(g,v_j) + d_G(v_j, v_i) + 1] \\
    &\quad - \sum_{1 \leq i < j \leq {2k+1}} \sum_{f \in E(H_i), g \in E(H_j)}[d_{G'}(f,v_1) + d_{G'}(g,v_1) + 1] \\
    &= \sum_{1 \leq i < j \leq 2k} \sum_{f \in E(H_i), g \in E(H_j)}d_G(v_i, v_j).
\end{split}
\end{equation*}
If $i \neq j$, $d_G(v_i, v_j) > 0$. So $S_3 > 0$.
\begin{equation*}
\begin{split}
S_4 &= \sum_{i=1}^{2k+1} \sum_{f \in E(H_i), g \in E(C_{2k+1})}[d_G(f,g) - d_{G'}(f,g)] \\
    &= \sum_{1 \leq i \leq {2k+1}, 1 \leq j \leq {2k+1}, g \in E(H_j)}d_G(e_i, g) - \sum_{1 \leq i \leq {2k+1}, 1 \leq j \leq {2k+1}, g \in E(H_j)}d_{G'}(e_i, g) \\
    &= \sum_{1 \leq i \leq {2k+1}, 1 \leq j \leq {2k+1}}m_jd_G(e_i, v_j) - \sum_{1 \leq i \leq {2k+1}, 1 \leq j \leq {2k+1}}m_jd_{G'}(e_i, v_1) \\
    &= \sum_{1 \leq j \leq {2k+1}}[m_j \sum_{i=1}^{2k+1} d_G(e_i, v_j)] - \sum_{1 \leq j \leq {2k+1}}[m_j \sum_{i=1}^{2k+1} d_{G'}(e_i, v_1)].
\end{split}
\end{equation*}
It can be checked that
\[
\sum_{i=1}^{2k+1}d_G(e_i, v_1) = [(0 + 1 + 2 + 3 + \cdots + k-1) \times 2 + k] = k^2.
\]
Then,
\begin{equation*}
\begin{split}
S_4 &= [m_1k^2 + m_2k^2 + m_3k^2 + \cdots + m_{2k}k^2 + m_{2k+1}k^2] \\
    &\quad - [m_1k^2 + m_2k^2 + m_3k^2 + \cdots + m_{2k}k^2 + m_{2k+1}k^2] \\
    &= 0.
\end{split}
\end{equation*}
Now we know $S_1 = S_2 = S_4 = 0$, so
\[
W_e(G) - W_e(G') = S_1 + S_2 + S_3 + S_4 = S_3 > 0
\]

That's the end of the proof.
\qed

\begin{lem}\label{lem-5}
Let $C_r = v_1v_2 \cdots v_rv_1$ $(r \geq 5)$ be an end-block of $G$ with $d_G(v_1) \geq 2$. Let $G' = G - v_{r-1}v_r - v_2v_3 + v_{r-1}v_1 + v_3v_1$; see Figure 4.
Let $m = |E(G_0)|$. Then we have $W_e(G') < W_e(G)$.

\end{lem}

\begin{figure}
\centering
    \includegraphics[height=6cm,width=10cm]{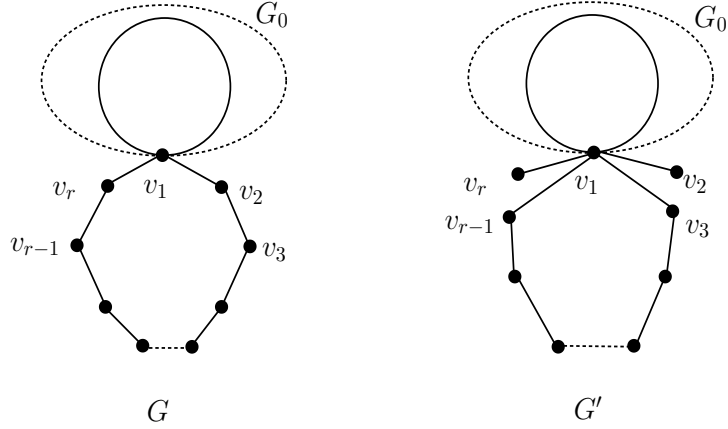}
    \caption{The graphs $G$ and $G'$ in Lemma 6}
    \label{fig4}
\end{figure}

\demo We calculate the edge-Wiener index of $G$ and $G'$, respectively.
We will deal with the problem with two cases according to the parity of $r$.

\vspace{2mm} {\bf  Case 1.}  $r$ is even.

By the definition of edge-Wiener, we can suppose that
\[
W_e(G) - W_e(G') = T_1 + T_2 + T_3 + T_4 + T_5 + T_6.
\]
Where,
\begin{equation*}
\begin{split}
T_1 &= \sum_{f,g \in E(G_0)} d_G(f,g) - \sum_{f,g \in E(G_0)}d_{G'}(f,g), \\
T_2 &= \sum_{f,g \in E(C_r)} d_G(f,g) - \sum_{f,g \in E(C_{r-2})}d_{G'}(f,g), \\
T_3 &= \sum_{f \in E(C_r), g \in E(G_0)}d_G(f,g) - \sum_{f \in E(C_{r-2}), g \in E(G_0)}d_{G'}(f,g), \\
T_4 &= 0 - 2\sum_{f \in E(C_{r-2})}d_{G'}(v_1v_2, f), \\
T_5 &= 0 - 2\sum_{f \in E(G_0)}d_{G'}(v_1v_2, f), \\
T_6 &= 0 - d_{G'}(v_1v_r,v_1v_2) = -1.
\end{split}
\end{equation*}
It is  not difficult to find that $T_1 = 0$ and the edge-Wiener index of cycle $C_r$ is equal to the Wiener index, so as to $W_e(C_r) = W(C_r)$. By the fact that the Wiener index of  a cycle equal to its edge-Wiener, then by Lemma~\ref{lem00}, one has that
\begin{equation*}
\begin{split}
T_2 &= \sum_{f,g \in E(C_r)}d_G(f,g) - \sum_{f,g \in E(C_{r-2})}d_{G'}(f,g) \\
    &= W_e(C_r) - W_e(C_{r-2}) \\
    &= \frac{3}{8}r^3 - \frac{3}{8}(r-3)^3. \\
\end{split}
\end{equation*}
By the symmetry of the cycle, we have
\begin{equation*}
\begin{split}
T_3 &= \sum_{f \in E(C_r), g \in E(G_0)}d_G(f,g) - \sum_{f \in E(C_{r-2}), g \in E(G_0)}d_{G'}(f,g) \\
    &= \sum_{f \in E(C_{r-2}), g \in E(G_0)}d_{G'}(f,g) + 2\sum_{f \in E(G_0)}d_{G'}(f,e_0) - \sum_{f \in E(C_{r-2}), g \in E(G_0)}d_{G'}(f,g) \\
    &= 2\sum_{f \in E(G_0)}d_{G'}(f,e_0)
\end{split}
\end{equation*}
with $e_0 = v_{\frac{r}{2}}v_{\frac{r}{2} + 1}$.
By calculation, we have
\begin{equation*}
\begin{split}
T_4 &= 0 - 2\sum_{f \in E(C_{r-2})}d_{G'}(v_1v_2,f) \\
    &= -2 \times 2 \times (1 + 2 + 3 + \cdots + \frac{r}{2} - 1) \\
    &= -4\sum_{i=1}^{\frac{r}{2} - 1} i,
\end{split}
\end{equation*}
\[
T_5 = -2\sum_{f \in E(G_0)}d_{G'}(v_1v_2,f), T_6 = -1.
\]
Thus,
\begin{equation*}
\begin{split}
W_e(G) - W_e(G') &= T_1 + T_2 + T_3 + T_4 + T_5 + T_6 \\
                 &= \frac{9}{4}r^2 - \frac{9}{2}r + 3 + 2\sum_{f \in E(G_0)}(\frac{r}{2} - 1) - \frac{1}{2}(r^2 - 2r) - 1 \\
                 &= \frac{7}{4}r^2 + (m - \frac{7}{2})r - 2m + 2  (r \geq 5).
\end{split}
\end{equation*}
We can get it easily that $W_e(G) - W_e(G') > 0$ for $r$ is even when $r \geq 5$.

\vspace{2mm} {\bf  Case 2.}  $r$ is odd.

By the definition of edge-Wiener, we can suppose that
\[
W_e(G) - W_e(G') = S_1 + S_2 + S_3 + S_4 + S_5 + S_6.
\]
Where,
\begin{equation*}
\begin{split}
S_1 &= \sum_{f,g \in E(G_0)}d_G(f,g) - \sum_{f,g \in E({G'}_0)}d_{G'}(f,g), \\
S_2 &= \sum_{f,g \in E(C_r)}d_G(f,g) - \sum_{f,g \in E(C_{r-2})}d_{G'}(f,g), \\
S_3 &= \sum_{f \in E(C_r), g \in E(G_0)}d_G(f,g) - \sum_{f \in E(C_{r-2}), g \in E(G_0)}d_{G'}(f,g), \\
S_4 &= 0 - 2\sum_{f \in E(C_{r-2})}d_{G'}(v_1v_2, f), \\
S_5 &= 0 - 2\sum_{f \in E(G_0)}d_{G'}(v_1v_2, f), \\
S_6 &= 0 - d_{G'}(v_1v_r,v_1v_2) = -1.
\end{split}
\end{equation*}
It can be checked that $S_1 = 0$, and $S_2 = \frac{1}{8}r(r^2 - 1) - \frac{1}{8}(r-2)[(r-2)^2 - 1]$,
\[
S_3 = 2\sum_{f \in E(G_0)}[d_{G'}(f,v_0) + (\frac{r-1}{2} - 1) + 1] + \sum_{f \in E(G_0)}(\frac{r-1}{2} - \frac{r-3}{2}),
\]
\[
S_4 = 0 - 2\sum_{f \in E(C_{r-2})}d_{G'}(v_1v_2, f) = -[2 \times (2\sum_{i=1}^\frac{r-3}{2} i + \frac{r-3}{2} + 1)].
\]
Then,
\begin{equation*}
\begin{split}
W_e(G) - W_e(G') &= S_1 + S_2 + S_3 + S_4 + S_5 + S_6 \\
                 &= \frac{1}{4}r^2 + (m - \frac{1}{2})r - 2m + \frac{3}{4}.
\end{split}
\end{equation*}
We can get it easily that $W_e(G) - W_e(G') > 0$ for $r$ is odd when $r \geq 5$.
\qed

\begin{lem}\label{lem-6}
Let $C_r = v_1v_2v_3v_4v_1$ be an end-block of $G$ with $d_G(v_1) \geq 2$. Let $G' = G - - v_1v_4 + v_2v_4$; see Figure 5.
Let $m = |E(G_0)|$. Then we have $W_e(G') < W_e(G)$.

\end{lem}

\begin{figure}
\centering
    \includegraphics[height=7cm,width=11cm]{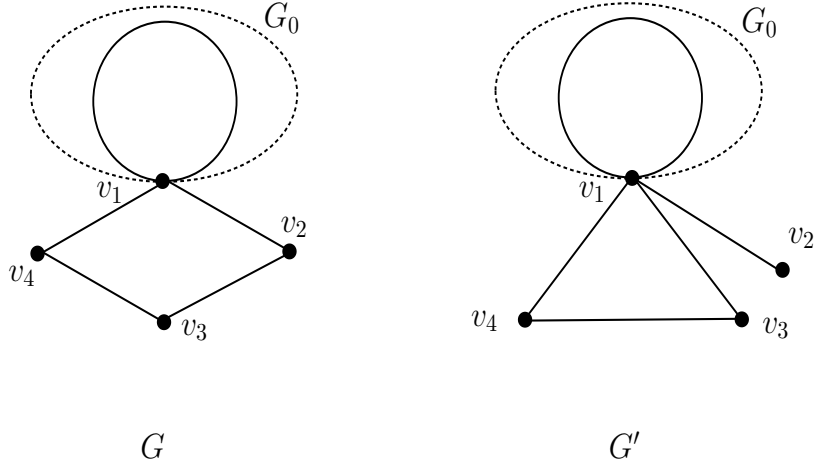}
    \caption{The graphs $G$ and $G'$ in Lemma 7}
    \label{fig5}
\end{figure}

\demo By definition, we have
\begin{equation*}
\begin{split}
W_e(G)-W_e(G')
      &= \sum_{f,g \in E(G_0)}[d_G(f,g) - d_{G'}(f,g)] + [\sum_{f,g \in E(C_4)}d_G(f,g) - \sum_{f,g \in E(C_3)} d_{G'}(f,g)] \\
      &\quad + [\sum_{f \in E(G_0), g \in E(C_4)}d_G(f,g) - \sum_{f \in E(G_0), g \in E(C_3)}d_{G'}(f,g)] \\
      &\quad - \sum_{f \in E(C_3)}d_{G'}(v_1v_2, f) - \sum_{f \in E(G_0)}d_{G'}(v_1v_2, f) \\
      &= \frac{4^3}{8} - \frac{3 \times (3^2-1)}{8} + \sum_{f \in E(G_0)}[d_{G}(f, v_1) + 1 + 1] - 4 - \sum_{f \in E(G_0)}[d_{G'}(f, v_1) + 1] \\
      &= 1 + \sum_{f \in E(G_0)}1 \\
      &= 1 + m > 0.
\end{split}
\end{equation*}

That's end of the proof.
\qed

Let $H_i$ be a connected graph with $x_i \in V(H_i)$ and suppose that $| V(H_i)| = n_i$, $| E(H_i) \mid = m_i$ for $i = 1,2,  \cdots ,l$ with $n_1 = \mbox{ max }\{n_1, n_2, \cdots ,n_l\}$, $m_1 = \mbox{ max }\{m_1, m_2, \cdots ,m_l\}$. Let $C_l = v_1v_2 \cdots v_lv_1$ be a cycle of length $l \geq 4$ and $D_l=C_l-v_1v_l+v_{l-2}v_l$. The graph $G$ is obtained from $C_l$ and $H_1,  \cdots , H_l$ by identifying $x_i$ with $v_i$ for $i = 1,2,  \cdots ,l$. The graph $G'$ is obtained from $G$ by deleting the edge $v_1v_l$ and adding an edge $v_{l-2}v_l$, i.e., $G' = G - v_1v_l + v_{l-2}v_l$.

\begin{lem}\label{lem-7}
Let $G$ and $G'$ be the above specified graphs; see Figure 6. Let $|E(H_i)| = m_i$ for $i=1, 2, \cdots, l$. Then we have

(1) If $l = 4$, then $W_e(G) \leq W_e(G')$ with equality if and only if $m_1 = m_2 = m_3 = m_4 = 1$.

(2) If $l \geq 5$, then $W_e(G) < W_e(G')$.

\end{lem}

\begin{figure}
\centering
    \includegraphics[height=5cm,width=6cm]{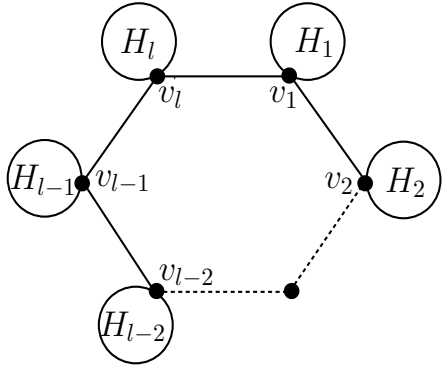}
    \caption{The graphs $G$ in Lemma 8}
    \label{fig6}
\end{figure}

\demo (1) $l = 4$. By definition and a direct calculation, one has that
\[
W_e(G') - W_e(G) = m_1m_4 - m_2m_4 + 2m_1 - m_2 - 1 = (m_1 - m_2)m_4 + 2m_1 -  m_2 - 1 \geq 0
\]
with equality if and only if $m_1 = m_2 = m_3 = m_4 = 1$.

(2) $l \geq 5$. We proceed by considering the parity of $l$ in the following manner.

\vspace{2mm} {\bf  Case 1.}  $l$ is even,  let $l=2k$ with $l \geq 3$.

Then according to the definition of $W_e(G)$ we obtain
\begin{equation*}
\begin{split}
W_e(G') - W_e(G) &= \sum_{f,g \in E(G')}d_{G'}(f,g) - \sum_{f,g \in E(G)}d_G(f,g) \\
                 &= \sum_{i=1}^{2k}\sum_{f,g \in E(H_i)}[d_{G'}(f,g) - d_G(f,g)] + \sum_{1 \leq i < j \leq 2k}\sum_{f \in E(H_i), g \in E(H_j)}[d_{G'}(f,g) - d_G(f,g)] \\
                 &\quad + \sum_{f,g \in E(D_l)}d_{G'}(f,g) - \sum_{f,g \in E(C_l)}d_G(f,g) \\
                 &\quad  + \sum_{1 \leq i \leq 2k, e_i \in  E(D_l), 1 \leq j \leq 2k, g \in  E(H_j)} d_{G'}(e_i,g) - \sum_{1 \leq i \leq 2k, e_i \in  E(C_l), 1 \leq j \leq 2k, g \in  E(H_j)} d_G(e_i,g).
\end{split}
\end{equation*}
We suppose $S_1$, $S_2$, $S_3$, $S_4$ in the following meaning
\begin{equation*}
\begin{split}
S_1 &= \sum_{i=1}^{2k}\sum_{f,g \in E(H_i)}[d_{G'}(f,g) - d_G(f,g)], \\
S_2 &= \sum_{1 \leq i < j \leq 2k}\sum_{f \in E(H_i), g \in E(H_j)}[d_{G'}(f,g) - d_G(f,g)], \\
S_3 &= \sum_{f,g \in E(D_l)}d_{G'}(f,g) - \sum_{f,g \in E(C_l)}d_G(f,g), \\
S_4 &= \sum_{1 \leq i \leq 2k, e_i \in  E(D_l), 1 \leq j \leq 2k, g \in  E(H_j)}d_{G'}(e_i,g) - \sum_{1 \leq i \leq 2k, e_i \in  E(C_l), 1 \leq j \leq 2k, g \in  E(H_j)}d_G(e_i,g). \\
\end{split}
\end{equation*}
It is not difficult to that $S_1 = 0$. Now let's calculate $S_2$,
\begin{equation*}
\begin{split}
S_2 &= \sum_{1 \leq i < j \leq 2k}\sum_{f \in E(H_i), g \in E(H_j)}[d_{G'}(f,g) - d_G(f,g)] \\
    &= \sum_{1 \leq i < j \leq 2k}m_im_j[d_{G'}(v_i,v_j) - d_G(v_i,v_j)] \\
    &= [\sum_{i=1}^{k-2}(\sum_{j=i+1}^{k+i} + \sum_{j=k+i+1}^{2k}) + \sum_{i=k-1}^{2k-2}(\sum_{j=i+1}^{2k-1} + \sum_{j=2k}) + \sum_{i=2k-1}\sum_{j=2k}] \\
    &\quad \times m_im_j[d_{G'}(v_i,v_j) - d_G(v_i,v_j)].
\end{split}
\end{equation*}
It can be checked that $d_{G'}(v_i,v_j) - d_G(v_i,v_j)=0$ holds for $(i,j) \in A \cup B \cup C$, where
\begin{equation*}
\begin{split}
A &= \{(i,j)|1 \leq i \leq k-2, i+1 \leq j \leq k+i\}, \\
B &= \{(i,j)|k-1 \leq i \leq 2k-2, i+1 \leq j \leq 2k-1\}, \\
C &= \{(i,j)|i = 2k-1, j = 2k\}.
\end{split}
\end{equation*}
Hence,
\begin{equation*}
\begin{split}
S_2 &= (\sum_{i=1}^{k-2}\sum_{j=k+i+1}^{2k} + \sum_{i=k-1}^{2k-2}\sum_{j=2k}) m_im_j[d_{G'}(v_i,v_j) - d_G(v_i,v_j)] \\
    &= \sum_{i=1}^{k-2}\sum_{j=k+i+1}^{2k-1} m_im_j(2j - 2i - 2k) + \sum_{i=1}^{k-2} m_im_{2k}(2k - 2i - 1) + m_{k-1}m_{2k} - \sum_{i=k}^{2k-2} m_im_{2k} \\
    &= \sum_{i=2}^{k-2}\sum_{j=k+i+1}^{2k-1} m_im_j(2j - 2i - 2k) + \sum_{i=2}^{k-2} m_im_{2k}(2k - 2i - 1) \\
    &\quad + \sum_{j=k+2}^{2k-1} m_1m_j(2j - 2k - 2) + m_1m_{2k}(2k - 3) + m_{k-1}m_{2k} - \sum_{i=k}^{2k-2} m_im_{2k}.
\end{split}
\end{equation*}
It obvious that $2j - 2i - 2k > 0$ for $k + i + 1 \leq j \leq 2k - 1$ and $2k - 2i - 1 > 0$ for $2 \leq i \leq k - 2$. This gives
\begin{equation*}
\begin{split}
S_2 &> \sum_{i=k+2}^{2k-1} m_1m_i(2i - 2k - 2) + m_1m_{2k}(2k - 3) + m_{k-1}m_{2k} - \sum_{i=k}^{2k-2} m_im_{2k} \\
    &= \sum_{i=k+2}^{2k-2} [m_1m_i(2i - 2k - 2) - m_im_{2k}] + m_1m_{2k-1}(2k - 4) + m_1m_{2k}(2k - 3) \\
    &+ m_{k-1}m_{2k} - m_km_{2k} - m_{k+1}m_{2k}.
\end{split}
\end{equation*}
Note that $m_1 =\mbox{ max } \{m_1, m_2, \cdots ,m_l\}$. Hence,
\begin{equation*}
\begin{split}
m_1m_i(2i - 2k - 2) - m_im_{2k} &\geq m_im_{2k}(2i - 2k - 2) - m_im_{2k} \\
                                &= m_im_{2k}(2i - 2k - 3) \\
                                &= m_im_{2k}[2(i - k - 2) + 1] > 0.
\end{split}
\end{equation*}
For $k + 2 \leq i \leq 2k - 2$, we have
\begin{equation*}
\begin{split}
S_2 &> m_1m_{2k-1}(2k - 4) + m_1m_{2k}(2k - 3) + m_{k-1}m_{2k} - m_km_{2k} - m_{k+1}m_{2k} \\
    &> m_1m_{2k}(2k - 3) - m_km_{2k} - m_{k+1}m_{2k} \\
    &= m_1m_{2k}(k - 1) - m_km_{2k} + m_1m_{2k}(k - 2) - m_{k+1}m_{2k} \\
    &\geq m_km_{2k}(k - 1) - m_km_{2k} + m_{k+1}m_{2k}(k - 2) - m_{k+1}m_{2k} \\
    &= m_km_{2k}(k - 2) + m_{k+1}m_{2k}(k - 3) > 0.
\end{split}
\end{equation*}
Thus, when $k \geq 3$, we have $S_2 > 0$.

Then we calculate $S_3$. By Lemma~\ref{lem00}, one has that if $l = 2k$, $W_e(C_l) = W(C_l) = W(C_{2k}) = \frac{1}{8}l^3 = \frac{1}{8}{(2k)}^3 = k^3$, $W_e(P_{l-2}) = W(P_{l-3}) = \frac{1}{6}(l -2)(l -3)(l -4)$. Suppose the sum distance of all edged in $D_l$ is  $W_e(D_l)$, $m_s = |E(C_3)| = 3$, $m_t = |E(P_{l-2})| = l - 3$. Then by Lemma~\ref{lem-1}, one has that
\[
W_e(D_l) = W_e(C_3) + W_e(P_{l-2}) + m_s\sum_{f \in E(P_{l-2})}d_{D_l}(f,u) + m_t\sum_{g \in E(C_3)}d_{D_l}(g,u) + m_sm_t
\]
so we have,
\begin{equation*}
\begin{split}
S_3 &= \sum_{f,g \in E(D_l)}d_{G'}(f,g) - \sum_{f,g \in E(C_l)}d_G(f,g) \\
    &= W_e(D_l) - W_e(C_l) \\
    &= \frac{4}{3}k^3 - \frac{13}{3}k + 5 - k^3 \\
    &= \frac{1}{3}k^3 - \frac{13}{3}k + 5 > 0.
\end{split}
\end{equation*}
$S_4$ maybe a little intractable,
\begin{equation*}
\begin{split}
S_4 &= \sum_{1 \leq i \leq 2k, e_i \in  E(D_l), 1 \leq j \leq 2k, g \in  E(H_j)}d_{G'}(e_i,g) - \sum_{1 \leq i \leq 2k, e_i \in  E(C_l), 1 \leq j \leq 2k, g \in  E(H_j)}d_G(e_i,g) \\
    &= \sum_{1 \leq j \leq 2k} [m_j\sum_{1 \leq i \leq 2k, e_i \in  E(D_l)}d_{G'}(e_i,v_j)] - \sum_{1 \leq j \leq 2k} [m_j\sum_{1 \leq i \leq 2k, e_i \in  E(C_l)}d_G(e_i,v_j)] \\
    &= \sum_{1 \leq j \leq 2k} [m_j\sum_{1 \leq i \leq 2k, e_i \in  E(D_l)}d_{G'}(e_i,v_j) - m_j\sum_{1 \leq i \leq 2k, e_i \in  E(C_l)}d_G(e_i,v_j)] \\
    &= \sum_{i=1}^{k-1}m_i[i^2 - (2k + 2)i + k^2 + 2k - 1] - m_k - 2m_{k+1} - m_{k+2} \\
    &\quad + \sum_{i=k+3}^{2k-2}m_i[i^2 - (2k + 2)i + k^2 + 2k - 1] + m_{2k-1}(k - 2)^2 + m_{2k}(k - 2)^2 \\
    &= m_1(k^2 - 2) + \sum_{i=2}^{k-1}m_i[i^2 - (2k + 2)i + k^2 + 2k - 1] - m_k - 2m_{k+1} - m_{k+2} \\
    &\quad + \sum_{i=k+3}^{2k-2}m_i[i^2 - (2k + 2)i + k^2 + 2k - 1] + m_{2k-1}(k - 2)^2 + m_{2k}(k - 2)^2.
\end{split}
\end{equation*}
When $k \geq 3$, we can calculate that $i^2 - (2k + 2)i + k^2 + 2k - 1 > 0$ when $1 \leq i \leq k - 1$, and $k + 3 \leq i \leq 2k - 2$, and $\mbox{ min } \{ k^2 - 2 \} = 7$, so $m_1(k^2 - 2) - m_k - 2m_{k+1} - m_{k+2} \geq 7m_1 - m_k - 2m_{k+1} - m_{k+2}$ with $m_1 =  \mbox{ max }\{m_1, m_2, \cdots ,m_l\}$, so $7m_1 - m_k - 2m_{k+1} - m_{k+2} > 0$, then $S_4 > 0$.

Then,
\[
W_e(G') - W_e(G) = S_1 + S_2 + S_3 + S_4 > 0.
\]
Thus, $W_e(G) < W_e(G')$ for $l = 2k$ with $k \geq 3$.

\vspace{2mm} {\bf  Case 2.}  $l$ is odd,  let $l=2k+1$ with $l \geq 2$.

Then according to the definition of edge-Wiener index, one has that
\begin{equation*}
\begin{split}
W_e(G') - W_e(G) &= \sum_{f,g \in E(G')}d_{G'}(f,g) - \sum_{f,g \in E(G)}d_G(f,g) \\
                 &= \sum_{i=1}^{2k+1}\sum_{f,g \in E(H_i)}[d_{G'}(f,g) - d_G(f,g)] \\
                 &\quad + \sum_{1 \leq i < j \leq 2k+1}\sum_{f \in E(H_i), g \in E(H_j)}[d_{G'}(f,g) - d_G(f,g)]\\
                 &\quad + \sum_{f,g \in E(D_l)}d_{G'}(f,g) - \sum_{f,g \in E(C_l)}d_G(f,g) \\
                 &\quad  + \sum_{1 \leq i \leq 2k+1, e_i \in  E(D_l), 1 \leq j \leq 2k+1, g \in  E(H_j)} d_{G'}(e_i,g) \\
                 &\quad - \sum_{1 \leq i \leq 2k+1, e_i \in  E(C_l), 1 \leq j \leq 2k+1, g \in  E(H_j)} d_G(e_i,g).
\end{split}
\end{equation*}
We suppose $T_1$, $T_2$, $T_3$, $T_4$ in the following meaning
\begin{equation*}
\begin{split}
T_1 &= \sum_{i=1}^{2k+1}\sum_{f,g \in E(H_i)}[d_{G'}(f,g) - d_G(f,g)], \\
T_2 &= \sum_{1 \leq i < j \leq 2k+1}\sum_{f \in E(H_i), g \in E(H_j)}[d_{G'}(f,g) - d_G(f,g)], \\
T_3 &= \sum_{f,g \in E(D_l)}d_{G'}(f,g) - \sum_{f,g \in E(C_l)}d_G(f,g), \\
T_4 &= \sum_{1 \leq i \leq 2k+1, e_i \in  E(D_l), 1 \leq j \leq 2k+1, g \in  E(H_j)}d_{G'}(e_i,g) - \sum_{1 \leq i \leq 2k+1, e_i \in  E(C_l), 1 \leq j \leq 2k+1, g \in  E(H_j)}d_G(e_i,g). \\
\end{split}
\end{equation*}
It is  not difficult to find that $T_1 = 0$. Now let's calculate $T_2$,
\begin{equation*}
\begin{split}
T_2 &= \sum_{1 \leq i < j \leq 2k+1}\sum_{f \in E(H_i), g \in E(H_j)}[d_{G'}(f,g) - d_G(f,g)] \\
    &= \sum_{1 \leq i < j \leq 2k+1}m_im_j[d_{G'}(v_i,v_j) - d_G(v_i,v_j)] \\
    &= [\sum_{i=1}^{k-1}(\sum_{j=i+1}^{k+i} + \sum_{j=k+i+1}^{2k+1}) + \sum_{i=k}^{2k-1}(\sum_{j=i+1}^{2k} + \sum_{j=2k+1}) + \sum_{i=2k}\sum_{j=2k+1}] \\
    &\quad \times m_im_j[d_{G'}(v_i,v_j) - d_G(v_i,v_j)].
\end{split}
\end{equation*}
It can be checked that $d_{G'}(v_i,v_j) - d_G(v_i,v_j)=0$ holds for $(i,j) \in {A'} \cup {B'} \cup {C'}$, where
\begin{equation*}
\begin{split}
A' &= \{(i,j)|1 \leq i \leq k-1, i+1 \leq j \leq k+i\} \\
B' &= \{(i,j)|k \leq i \leq 2k-1, i+1 \leq j \leq 2k\} \\
C' &= \{(i,j)|i = 2k, j = 2k+1\}
\end{split}
\end{equation*}
Hence,
\begin{equation*}
\begin{split}
T_2 &= (\sum_{i=1}^{k-1}\sum_{j=k+i+1}^{2k+1} + \sum_{i=k}^{2k-1}\sum_{j=2k+1}) m_im_j[d_{G'}(v_i,v_j) - d_G(v_i,v_j)] \\
    &= \sum_{i=1}^{k-1}\sum_{j=k+i+1}^{2k} m_im_j(2j - 2i - 2k - 1) + \sum_{i=1}^{k-1} m_im_{2k+1}(2k - 2i) - \sum_{i=k+1}^{2k-1} m_im_{2k+1} \\
    &= \sum_{i=2}^{k-1}\sum_{j=k+i+1}^{2k} m_im_j(2j - 2i - 2k - 1) + \sum_{i=2}^{k-1} m_im_{2k+1}(2k - 2i) \\
    &\quad + \sum_{j=k+2}^{2k} m_1m_j(2j - 2k - 3) + m_1m_{2k+1}(2k - 2) - \sum_{i=k+1}^{2k-1} m_im_{2k+1}.
\end{split}
\end{equation*}
It obvious that $2j - 2i - 2k - 1 > 0$ for $k + i + 1 \leq j \leq 2k$ and $2k - 2i > 0$ for $2 \leq i \leq k - 1$. This gives
\begin{equation*}
\begin{split}
T_2 &> \sum_{i=k+2}^{2k}m_1m_i(2i - 2k - 3) + m_1m_{2k+1}(2k - 2) - \sum_{i=k+1}^{2k-1}m_im_{2k+1} \\
    &= \sum_{i=k+2}^{2k-1} [m_1m_i(2i - 2k - 3) - m_im_{2k+1}] \\
    &+ m_1m_{2k}(2k - 3) + m_1m_{2k+1}(2k - 2) - m_{k+1}m_{2k+1}.
\end{split}
\end{equation*}
Note that $m_1 = \mbox{ max }\{m_1, m_2, \cdots ,m_l\}$. Hence,
\begin{equation*}
\begin{split}
m_1m_i(2i - 2k - 3) - m_im_{2k+1} &\geq m_im_{2k+1}(2i - 2k - 3) - m_im_{2k+1} \\
                                  &= m_im_{2k+1}(2i - 2k - 4) \\
                                  &= 2m_im_{2k+1}(i - k - 2) \geq 0
\end{split}
\end{equation*}
for $k + 2 \leq i \leq 2k - 1$. Therefore, we obtain
\begin{equation*}
\begin{split}
T_2 &> m_1m_{2k}(2k - 3) + m_1m_{2k+1}(2k - 2) + m_{k+1}m_{2k+1} \\
    &\geq m_1m_{2k}(2k - 3) + m_{k+1}m_{2k+1}(2k - 2) - m_{k+1}m_{2k+1} \\
    &= m_1m_{2k}(2k - 3) + m_{k+1}m_{2k+1}(2k - 3).
\end{split}
\end{equation*}
Then, we have $T_2 > 0$ when $k \geq 2$.

Then we calculate $T_3$.  Then by Lemma~\ref{lem00}, one has that, when $l = 2k + 1$, $W_e(C_l) = W(C_l) = W(C_{2k+1}) = \frac{1}{8}(l^2 - 1) = \frac{1}{8}k(k + 1)(2k + 1)$, $W_e(P_{l-2}) = W(P_{l-3}) = \frac{1}{6}(l -2)(l -3)(l -4)$. Suppose the sum distance of all edged in $D_l$ is  $W_e(D_l)$, $m_s = |E(C_3)| = 3$, $m_t = |E(P_{l-2})| = l - 3$. Then by Lemma~\ref{lem-1}, one has that
\[
W_e(D_l) = W_e(C_3) + W_e(P_{l-2}) + m_s\sum_{f \in E(P_{l-2})}d_{D_l}(f,u) + m_t\sum_{g \in E(C_3)}d_{D_l}(g,u) + m_sm_t
\]
so we have,
\begin{equation*}
\begin{split}
T_3 &= \sum_{f,g \in E(D_l)}d_{G'}(f,g) - \sum_{f,g \in E(C_l)}d_G(f,g) \\
    &= W_e(D_l) - W_e(C_l) \\
    &= \frac{4}{3}k^3 + 2k^2 - \frac{10}{3}k + 3 - \frac{1}{2}k(k + 1)(2k + 1) \\
    &= \frac{1}{3}k^3 + \frac{1}{2}k^2 - \frac{23}{6}k + 3 > 0.
\end{split}
\end{equation*}
$T_4$ maybe a little intractable,
\begin{equation*}
\begin{split}
T_4 &= \sum_{1 \leq i \leq 2k+1, e_i \in  E(D_l), 1 \leq j \leq 2k+1, g \in  E(H_j)}d_{G'}(e_i,g) - \sum_{1 \leq i \leq 2k+1, e_i \in  E(C_l), 1 \leq j \leq 2k+1, g \in  E(H_j)}d_G(e_i,g) \\
    &= \sum_{1 \leq j \leq 2k+1}[m_j\sum_{1 \leq i \leq 2k+1, e_i \in  E(D_l)}d_{G'}(e_i,v_j)] - \sum_{1 \leq j \leq 2k+1} [m_j\sum_{1 \leq i \leq 2k+1,e_i \in  E(C_l)}d_G(e_i,v_j)] \\
    &= \sum_{1 \leq j \leq 2k+1}[m_j\sum_{1 \leq i \leq 2k+1, e_i \in  E(D_l)}d_{G'}(e_i,v_j) - m_j\sum_{1 \leq i \leq 2k+1, e_i \in  E(C_l)}d_G(e_i,v_j)] \\
    &= \sum_{i=1}^{k-1}m_i[i^2 - (2k + 3)i + k^2 + 3k] - 2m_{k+1} - 2m_{k+2} \\
    &\quad+ \sum_{i=k+3}^{2k-1}m_i[i^2 - (2k + 3)i + k^2 + 3k] + m_{2k}(k - 1)(k - 2) + m_{2k+1}(k - 1)(k - 2) \\
    &= m_1(k^2 + k - 2) + \sum_{i=2}^{k-1}m_i[i^2 - (2k + 3)i + k^2 + 3k] - 2m_{k+1} - 2m_{k+2} \\
    &+ \sum_{i=k+3}^{2k-1}m_i[i^2 - (2k + 3)i + k^2 + 3k] + m_{2k}(k - 1)(k - 2) + m_{2k+1}(k - 1)(k - 2).
\end{split}
\end{equation*}
When $k \geq 2$, we can calculate that $i^2 - (2k + 3)i + k^2 + 3k \geq 0$ when $i \leq i \leq k - 1$ and $k + 3 \leq i \leq 2k - 1$, and $\mbox{ min } \{ k^2 + k - 2 \} = 4$, so $m_1(k^2 + k - 2) - 2m_{k+1} - 2m_{k+2} \geq 4m_1 - 2m_{k+1} - m_{k+2}$ with $m_1 = \mbox{ max }\{m_1, m_2, \cdots ,m_l\}$, so $4m_1 - 2m_{k+1} - m_{k+2} \geq 0$, then $T_4 > 0$. Thus
\[
W_e(G') - W_e(G) = T_1 + T_2 + T_3 + T_4 > 0.
\]
Thus, $W_e(G) < W_e(G')$ for $l = 2k + 1$ with $k \geq 2$.

That's end of the proof.
\qed

Let $G_i$ be a connected $n_i$-vertex graph satisfying $| V(G_i) | = n_i$, $| E(G_i) | = m_i$, and $n_i \geq 2$, $i = 1,2,3$, $m_0 = |E(G_0)|$. Denote the farthest vertex from vertex $u_1$ in $G_1$ by $v$ and the farthest vertex from $u_2$ in $G_2$ by $w$. Let $G_0$ be the graph obtained from $G_1$, $G_2$ by identifying $u_1$ with $u_2$ as a new vertex $u$. Then $G$ and $G'$ are the graphs obtained from the graphs $G_0$ and $G_3$ by identifying $u_3$ with $u$ and $v$, respectively. The graphs $G$ and $G'$ are depicted in Figure 7.

\begin{lem}\label{lem-8}
Let $G$, $G'$ be the above specified graphs; see Figure 7. Then $W_e(G) \leq W_e(G')$ with equality if and only if $m_2 \geq m_1$.

\end{lem}

\begin{figure}
\centering
    \includegraphics[height=4cm,width=10cm]{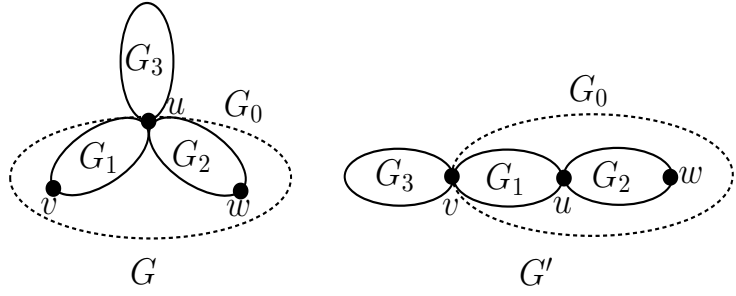}
    \caption{The graph $G$ and $G'$ in Lemma 9}
    \label{fig7}
\end{figure}

\demo We calculate the edge-Wiener index of $G$ and $G'$, respectively. By Lemma~\ref{lem-1},
\[
W_e(G) = W_e(G_0) + W_e(G_3) + m_3\sum_{f \in E(G_0)}d_G(f,u) + m_0\sum_{g \in E(G_3)}d_G(g,u) + m_0m_3
\]
\[
W_e(G') = W_e(G_0) + W_e(G_3) + m_3\sum_{f \in E(G_0)}d_{G'}(f,v) + m_0\sum_{g \in E(G_3)}d_{G'}(g,v) + m_0m_3
\]
\begin{equation*}
\begin{split}
W_e(G') - W_e(G) &= m_3[\sum_{f \in E(G_0)}d_{G'}(f,v) - \sum_{f \in E(G_0)}d_G(f,u)] \\
                 &\quad+ m_0[\sum_{g \in E(G_3)}d_{G'}(g,v) - \sum_{g \in E(G_3)}d_G(g,u)] \\
                 &= m_3[\sum_{f \in E(G_1)}d_{G'}(f,v) + \sum_{f \in E(G_2)}d_{G'}(f,v)] \\
                 &\quad- m_3[\sum_{f \in E(G_1)}d_G(f,u) + \sum_{f \in E(G_2)}d_G(f,u)] \\
                 &= m_3\{\sum_{f \in E(G_1)}[d_{G'}(f,v) - d_{G}(f,u)] + m_2d_{G_1}(u,v)\} \\
                 &\geq m_3[-m_1d_{G_1}(u,v) + m_2d_{G_1}(u,v)] \\
                 &= m_3(m_2 - m_1)d_{G_1}(u,v)
\end{split}
\end{equation*}
and $m_2 \geq m_1$, so $W_e(G') - W_e(G) \geq 0$.

That's the end of the proof.
\qed

Let $G_i$ be a connected $n_i$-vertex graph satisfying $| V(G_i) | = n_i$, $| E(G_i) | = m_i$, and $n_i \geq 2$, for $i = 1,2,3$, $m_0 = |E(G_0)|$. For any vertex $x_i \in V(G_i)$, denote the farthest vertex from vertex $x_2$ in $G_2$ by $u_2$. The graph $G_0$ is obtained from $G_2$, $G_3$ and $C_3 = v_1v_2v_3v_1$ by identifying $x_i$ with $v_i$ for $i = 2,3$. Then the graph $G$ and $G'$ are obtained from $G_0$ and $G_1$ by identifying $x_1$ with $v_1$ and $u_2$, respectively. The graphs $G$ and $G'$ are depicted in Figure 8.

\begin{lem}\label{lem-9}
Let $G$, $G'$ be the above specified graphs; see Figure 8. Then $W_e(G) < W_e(G')$ with $m_3 \geq m_2$.
\end{lem}

\begin{figure}
\centering
    \includegraphics[height=4cm,width=10cm]{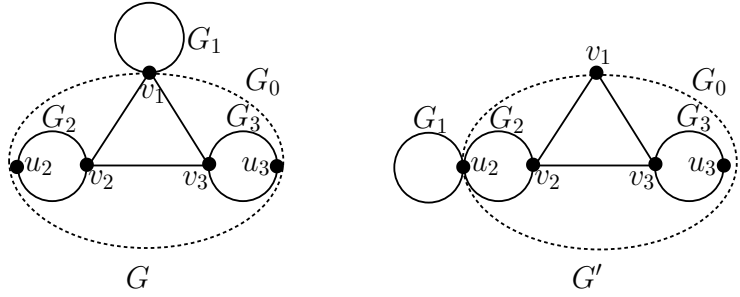}
    \caption{The graph $G$ and $G'$ in Lemma 10}
    \label{fig8}
\end{figure}

\demo  We calculate the edge-Wiener index of $G$ and $G'$, respectively. By Lemma~\ref{lem-1},
\[
W_e(G) = W_e(G_1) + W_e(G_0) + m_1\sum_{f \in E(G_0)}d_G(f,v_1) + m_0\sum_{g \in E(G_1)}d_G(g,v_1) + m_0m_1.
\]
\[
W_e(G') = W_e(G_1) + W_e(G_0) + m_1\sum_{f \in E(G_0)}d_{G'}(f,u_2) + m_0\sum_{g \in E(G_1)}d_{G'}(g,u_2) + m_0m_1.
\]
Then,
\begin{equation*}
\begin{split}
W_e(G') - W_e(G) &= m_1[\sum_{f \in E(G_0)}d_{G'}(f,u_2) - \sum_{f \in E(G_0)}d_G(f,v_1)] \\
                 &= m_1[\sum_{f \in E(G_2)}d_{G'}(f,u_2) + \sum_{f \in E(G_3)}d_{G'}(f,u_2) + d_{G'}(v_1v_2,u_2) + d_{G'}(v_1v_3,u_2) + d_{G'}(v_2v_3,u_2)] \\
                 &\quad- m_1[\sum_{f \in E(G_2)}d_G(f,v_1) + \sum_{f \in E(G_3)}d_G(f,v_1) + d_G(v_1v_2,v_1) + d_G(v_1v_3,v_1) + d_G(v_2v_3,v_1)] \\
                 &= m_1\{\sum_{f \in E(G_2)}[d_{G'}(f,u_2) - d_{G}(f,v_1)] + \sum_{f \in E(G_3)}[d_{G'}(f,u_2) - d_{G}(f,v_1)]\} \\
                 &\quad+ m_1[d_{G'}(u_2,v_2) + 1 + d_{G'}(u_2,v_2) + d_{G'}(u_2,v_2) - 1] \\
                 &= m_1\{\sum_{f \in E(G_2)} [d_{G'}(f,u_2) - 1 - d_{G}(f,v_2)] + m_3d_{G'}(u_2,v_2) + 3d_{G'}(u_2,v_2)\} \\
                 &\geq m_1\{-m_2[d_{G'}(v_2,u_2) - 1] + m_3d_{G'}(u_2,v_2) + 3d_{G'}(u_2,v_2)\} \\
                 &= m_1[(m_3 - m_2)d_{G'}(u_2,v_2) + m_2 + 3d_{G'}(u_2,v_2)].
\end{split}
\end{equation*}
By $m_3 \geq m_2$, so $W_e(G') - W_e(G) >0$.

That's the end of the proof.
\qed

For $v \in V(G)$, we define $$D_v(G)=\sum_{f\in E(G)}d_G(v, f).$$

Consider any two connected graphs $G_1$, $G_2$, where $u_1 \in V(G_1)$ and $| V(G_1) | = n_i \geq 2$, $| E(G_i) | = m_i$, $i = 1,2$, $m_0 = |E(G_0)| = |E({G'}_0)|$, and $G_2$ is not a path. Let $u$, $v$ be the end-vertices of a longest path of $G_2$. Clearly, $n_2 = | V(G_2) | \geq 3$. The graph $G$ is obtained from $G_1$, $G_2$ and $P_s = p_1p_2 \cdots p_s$ by identifying $u_1$ with $u$ and $v$ with $p_1$, respectively. The graph $G'$ is obtained from $G_1$, $G_2$ and $P_s = p_1p_2 \cdots p_s$ in a similar manner, which are depicted in Figure 9.

\begin{lem}\label{lem-10}
Let $G$, $G'$ be the above specified graphs; see Figure 9. Then $W_e(G) < W_e(G')$ with $D_v(G_2) \geq D_u(G_2)$.
\end{lem}

\begin{figure}
\centering
    \includegraphics[height=7cm,width=13cm]{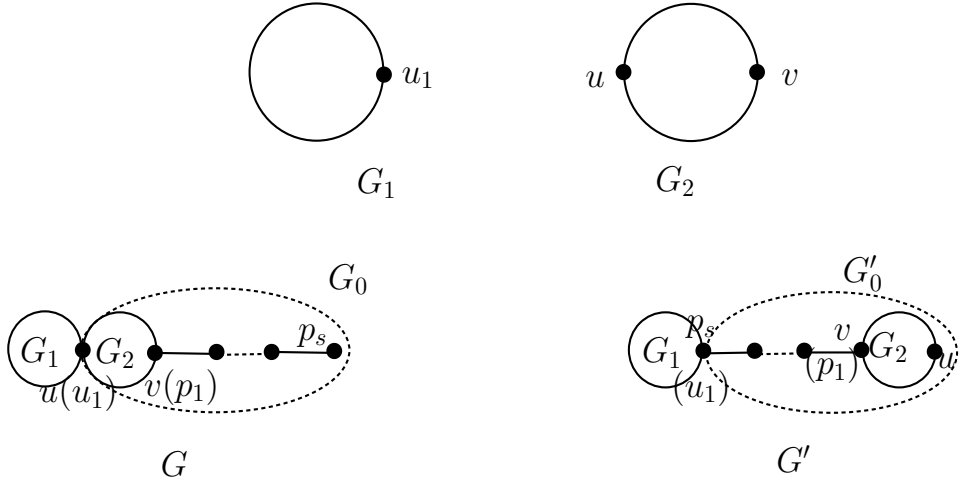}
    \caption{The graph $G$ and $G'$ in Lemma 11}
    \label{fig9}
\end{figure}

\demo  We calculate the edge-Wiener index of $G$ and $G'$, respectively. By Lemma~\ref{lem-1},
\[
W_e(G) = W_e(G_1) + W_e(G_0) + m_1\sum_{f \in E(G_0)}d_G(f,u) + m_0\sum_{g \in E(G_1)}d_G(g,u) + m_0m_1.
\]
\[
W_e(G') = W_e(G_1) + W_e({G'}_0) + m_1\sum_{f \in E({G'}_0)}d_{G'}(f,u_1) + m_0\sum_{g \in E(G_1)}d_{G'}(g,u_1) + m_0m_1.
\]
According to the structure of the graph, $m_0\sum_{g \in E(G_1)} d_G(g,u) = m_0\sum_{g \in E(G_1)} d_{G'}(g,u_1)$. So
\[
W_e(G') - W_e(G) = m_1\sum_{f \in E({G'}_0)} d_{G'}(f,u_1) - m_1\sum_{f \in E(G_0)} d_G(f,u).
\]
Let the length of the path $P_s = p_1p_2 \cdots p_s$ be $s$, $s = |E(P_s)|$. So,
\begin{equation*}
\begin{split}
W_e(G') - W_e(G) &= m_1\sum_{f \in E({G'}_0)}d_{G'}(f,u_1) - m_1\sum_{f \in E(G_0)}d_G(f,u) \\
                 &= m_1\{\sum_{f \in E(P_s)}d_{G'}(f,u_1) + \sum_{f \in E(G_2)}[d_{G'}(f,v) + s - 1]\} \\
                 &\quad- m_1\{\sum_{f \in E(G_2)}d_G(f,u) + \sum_{f \in E(P_s)}[d_G(f,v) + d_{G}(u,v)]\} \\
                 &= m_1[\sum_{f \in E(G_2)}d_{G'}(f,v) + m_2(s - 1) - \sum_{f \in E(G_2)}d_G(f,u) - (s - 1)d_{G}(u,v)] \\
                 &= m_1\{\sum_{f \in E(G_2)}[d_{G'}(f,v) - d_{G}(f,u)] + (s - 1)[m_2 - d_{G}(u,v)]\}.
\end{split}
\end{equation*}
Consider $G_2$ is not a path, so $m_2 - d(u,v) > 0$, and $D_v(G_2) \geq D_u(G_2)$. Then $\sum_{f \in E(G_2)}[d(f,v) - d(f,u)] \geq 0$, which we obtain $W_e(G') - W_e(G) > 0$.

That's the end of the proof.
\qed

Let $G_1$, $G_2$ be two graphs with $u \in V(G_1)$, $V \in V(G_2)$, satisfying $| V(G_1) | = n_i \geq 2$, $| E(G_i) | = m_i$, $i = 1,2$, $m_0 = |E(G_0)| = |E({G'}_0)|$. We assume that $m_2 \geq m_1 \geq 1$. Denote a path with length $s \geq 4$ by $P_s = p_1p_2 \cdots p_s$. Then put $H = P_s + p_{s-2}p_s$ and $H' = H - p_{s-2}p_s - p_{s-1}p_s + p_1p_s + p_2p_s$. The graph $G$ is obtained from $G_1$, $G_2$ and $H$ by identifying $u$ with $p_1$ (resp. $v$ with $p_{s-1}$), whereas $G'$ is obtained from $G_1$, $G_2$ and $H'$ by identifying $u$ with $p_1$ (resp. $v$ with $p_{s-1}$). The graph $H$, $H'$, $G$ and $G'$ are depicted in Figure 10.

\begin{lem}\label{lem-11}
Let $G$, $G'$ be the above specified graphs; see Figure 10. Then $W_e(G) \leq W_e(G')$ with equality if and only if $m_2 = m_1$.

\end{lem}

\begin{figure}
\centering
    \includegraphics[height=5cm,width=12cm]{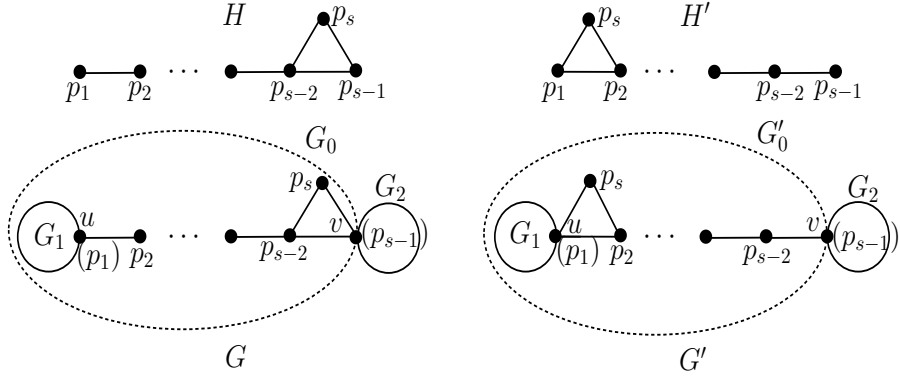}
    \caption{The graph $G$ and $G'$ in Lemma 12}
    \label{fig10}
\end{figure}

\demo  We calculate the edge-Wiener index of $G$ and $G'$, respectively. By Lemma~\ref{lem-1},
\[
W_e(G) = W_e(G_2) + W_e(G_0) + m_0\sum_{f \in E(G_2)}d_G(f,v) + m_2\sum_{g \in E(G_0)}d_G(g,v) + m_0m_2.
\]
\[
W_e(G') = W_e(G_2) + W_e({G'}_0) + m_0\sum_{f \in E(G_2)}d_{G'}(f,v) + m_2\sum_{g \in E({G'}_0)}d_{G'}(g,v) + m_0m_2.
\]
Subtracting the two type,
\[
W_e(G') - W_e(G) = W_e({G'}_0) - W_e(G_0) + m_2[\sum_{g \in E({G'}_0)} d_{G'}(g,v) - \sum_{g \in E(G_0)} d_G(g,v)].
\]
Now let's calculate $W_e({G'}_0) - W_e(G_0)$ and $m_2[\sum_{g \in E({G'}_0)} d_{G'}(g,v) - \sum_{g \in E(G_0)} d_G(g,v)]$ separately. By Lemma~\ref{lem-1},
\[
W_e({G'}_0) = W_e(G_1) + W_e(G_{H'}) + m_1\sum_{f \in E(H')}d_{G'_0}(f,u) +  m_{H'}\sum_{g \in E(G_1)}d_{G'_0}(g,u) + m_1m_{H'}.
\]
\[
W_e(G_0) = W_e(G_1) + W_e(G_H) + m_1\sum_{f \in E(H)}d_{G_0}(f,u) +  m_H\sum_{g \in E(G_1)}d_{G_0}(g,u) + m_1m_H.
\]
\begin{equation*}
\begin{split}
W_e({G'}_0) - W_e(G_0) &= m_1[\sum_{f \in E(H')}d_{G'_0}(f,u) - \sum_{f \in E(H)}d_{G_0}(f,u)] \\
                       &= m_1\{(1 + \sum_{i=0}^{s-1-2} i) - [\sum_{i=0}^{s-1-2} i + (s - 1 - 2) + (s - 1 - 1)]\} \\
                       &= m_1(6 - 2s).
\end{split}
\end{equation*}
then we calculate $m_2[\sum_{g \in E({G'}_0)}d_{G'}(g,v) - \sum_{g \in E(G_0)}d_G(g,v)]$. Suppose $M = m_2[\sum_{g \in E({G'}_0)}d_{G'}(g,v) - \sum_{g \in E(G_0)}d_G(g,v)]$, then
\begin{equation*}
\begin{split}
   M &= m_2[\sum_{g \in E(G_1)}d_{G'}(g,v) + \sum_{g \in E(H')}d_{G'}(g,v) - \sum_{g \in E(G_1)}d_G(g,v) -  \sum_{g \in E(H)}d_G(g,v)] \\
     &= m_2[\sum_{g \in E(H')}d_{G'}(g,v) - \sum_{g \in E(H)}d_G(g,v)] \\
     &= m_2[\sum_{i=0}^{s-1-2} i + (s - 1 - 2) + (s - 1 - 1) - \sum_{i=0}^{s-1-2} i - 1] \\
     &= m_2(2s - 6).
\end{split}
\end{equation*}
So we have
\[
W_e(G') - W_e(G) = m_1(6 - 2s) + m_2(2s - 6) = (m_2 - m_1)(2s - 6).
\]
If $m_2 > m_1$ and $s \geq 4$, $(m_2 - m_1)(2s - 6) \geq 0$, one has that $W_e(G') - W_e(G) > 0$. When $m_2 = m_1$, we obtain $W_e(G') = W_e(G)$.

That's the end of the proof.
\qed

\section{Cactus with minimum and maxmimum edge-Wiener indices in $\mathcal{G}_{n,t}$}
In this section, based on the results obtained in Section 2, we establish a sharp upper and a lower bounds on the sum of all distances of the graphs in $\mathcal{G}_{n,t}$.

\begin{theorem}\label{cp}
For any $G \in \mathcal{G}_{n,t}$, one has that
\[
W_e(G) \geq \frac{1}{2}n^2 + (2t - \frac{3}{2})n + 3t^2 - 7t + 1
\]
with equality if and only if $G \cong C_0(n,t)$; see Figure 11.
\end{theorem}

\begin{figure}
\centering
    \includegraphics[height=4cm,width=7cm]{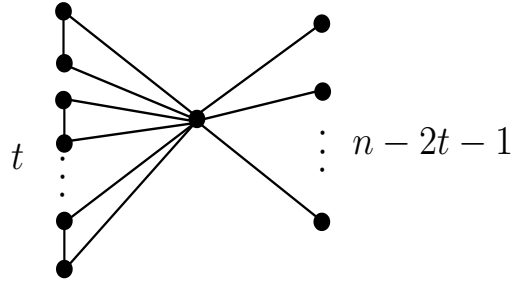}
    \caption{The graph $C_0(n,t)$}
    \label{fig11}
\end{figure}

\demo Supposed that $G$ is the graph that has minimum edge-Wiener index in $\mathcal{G}_{n,t}$. By Lemma~\ref{lem-2}, we have that all the cut edges of $G$ are pendent edges. By Lemmas~\ref{lem-3} and~\ref{lem-4}, we can see that all the cycles of $G$ are end-block. So we have that there exist $g_1, g_2, \cdots ,g_k$ such that $G \cong G(g_1, g_2, \cdots ,g_k)$. By Lemmas~\ref{lem-5} and~\ref{lem-6}, we have that $g_i = 3$ for $i = 1,2, \cdots ,k$. Hence, $G \cong C_0(n,t)$.
By simple calculation, we have that for any $G \in \mathcal{G}_{n,t}$,
\[
W_e(G) \geq \frac{1}{2}n^2 + (2t - \frac{3}{2})n + 3t^2 - 7t + 1
\]
with equality if and only if $G \cong C_0(n,t)$.
\qed

Let $H_1$ (resp. $H_2$) be a triangle-chain of length $i$ (resp. $j$) with an end
$u$ (resp. $v$). Joining the end $u$ of $H_1$ with the end $v$ of $H_2$ by a path $P_{n-2i-2j}$ yields the saw graph $Sw( i, j; n-2i-2j-1   )$.
The saw graph $Sw( \lfloor   \frac{t}{2}  \rfloor,  \lceil   \frac{t}{2} \rceil; n-2t-1   )$ is depicted in Figure 12.

\begin{theorem}\label{cp}
For any $G \in \mathcal{G}_{n,t}$ such that $n \geq 5$, and $t \geq 0$,one has that

(i)If $t = 2k$ $(k \geq 0)$, then
\[
W_e(G) \leq \frac{1}{6}n^3 - \frac{1}{2}n^2 + \frac{1}{3}n + kn^2 - 4k^2n + \frac{4}{3}k^3 + 8k^2 - \frac{10}{3}k
\]
with equality if and only if $G \cong Sw(k,k;n - 4k - 1)$.

(ii)If $t = 2k + 1$ $(k \geq 0)$, then
\[
W_e(G) \leq \frac{1}{6}n^3 - \frac{13}{6}n + kn^2 - 4k^2n - 4kn - \frac{4}{3}k^3 + 10k^2 + \frac{35}{3}k + 5
\]
with equality if and only if $G \cong Sw(k,k+1;n - 4k - 3)$.
\end{theorem}

\begin{figure}
\centering
    \includegraphics[height=2cm,width=10cm]{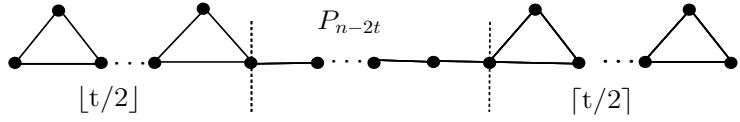}
    \caption{The graph $Sw(\lfloor {t/2} \rfloor,\lceil {t/2} \rceil;n - 2t - 1)$}
    \label{fig12}
\end{figure}

\demo Choose $G$ in $G \in \mathcal{G}_{n,t}$ with $n \geq 5$, and $t \geq 0$ such that its sum of distances is as large as possible. We proceed by considering whether $G$ contains cycles or not.

If $G$ does not contain a cycle, i.e., $t=0$, then $G$ is a tree. Suppose that $G$ is not the path. Then apply Lemma~\ref{lem-8} repeatedly
until $G$ contains exactly two pendant vertices. Thus, we get
$$W_e(G) < W_e(P_n)=W_e(Sw(0, 0; n-1))=\frac{1}{6}n^{3}-\frac{1}{2}n^{2}+\frac{1}{3}n.$$

If $G$ contains cycle(s), i.e., $t\geq 1$, then $G$ is in $ \mathcal{G}_{n,t}$ with $n\geq 5$ and $t \geq 1$.

\begin{clm}\label{cp}
Each of the cycles contained in $G$ is of length 3.
\end{clm}

\noindent\textbf{Proof of Claim 1.}
By Lemma~\ref{lem-7}, $G$ does not contain a cycle $C_l$ with $l \geq 5$; otherwise there exists another graph $H$ in $\mathcal{G}_{n,t}$ such that $W_e(G) < W_e(H)$, a contradiction to the choice of $G$. In fact, if $G$ contains a cycle $C_4 = v_1v_2v_3v_4v_1$, then we denote the components of $G - \{v_1v_2, v_2v_3, v_3v_4, v_1v_4\}$ containing $v_1$, $v_2$, $v_3$, and $v_4$, respectively, by $H_1$, $H_2$, $H_3$, and $H_4$ with $| E(H_1) | \geq | E(H_i) |$, for $i = 2, 3, 4$. If $| E(H_1) | \neq | E(H_i) |$, then according to Lemma~\ref{lem-7}, we can get a new graph $G' \in \mathcal{G}_{n,t}$ such that $W_e(G) < W_e(G')$, a contradiction to the choice of $G$. Otherwise, based on the graph transformation in Lemma~\ref{lem-7}, we have $W_e(G) = W_e(G')$. Therefore, each of the cycles contained in $G$ must be of length 3.
\qed

\begin{clm}\label{cp}
The graph $G$ is a chain cactus.
\end{clm}

\noindent\textbf{Proof of Claim 2.}
If there is some block in $G$ such that this block contains three cut vertices (based on Claim 1 every in $G$ is a triangle), then by Lemma~\ref{lem-9}, there exists a new graph $G'$ in $\mathcal{G}_{n,t}$ satisfying $W_e(G) < W_e(G')$, contradiction. Hence, each block in $G$ has at most two cut vertices. If there is some cut vertex in $G$ such that this vertex is shared by at least three blocks, then by Lemma~\ref{lem-8}, there is a graph $G'$ in $\mathcal{G}_{n,t}$, contradiction. Thus, each cut vertex in $G$ is shared by exactly two blocks. That is to say, $G$ is a chain cactus.
\qed

\begin{clm}\label{cp}
If $t = 1$, then the graph $G$ has exactly one pendant path. If $t \geq 2$, then $G$ has no pendant paths.
\end{clm}

\noindent\textbf{Proof of Claim 3.}
If $G$ contains a single cycle, then by Lemma~\ref{lem-8} we may suppose that there are two pendant paths attached to two different vertices of this cycle. By Lemma ~\ref{lem-10}, there exists a new graph $G'$ in $\mathcal{G}_{n,t}$ such that $W_e(G) < W_e(G')$ where $G'$ has only one pendant path, a contradiction. If $G$ contains at least two cycles, then we may suppose that $G$ contains pendant path. Note that since $G$ is a cactus chain (by Claim 2), the pendant path must be attached at the leftmost (or rightmost) triangle of $G$. Using Lemma~\ref{lem-10} twice, we can obtain a graph $G' \in \mathcal{G}_{n,t}$ such that $W_e(G) < W_e(G')$, where $G'$ contains no pendant path, a contradiction to the choice of $G$.
\qed

Let $v$ be a vertex of the graph $G$. The degree of $v$, denote by $d(v)$, is the number of first neighbors of this vertex. Hoffman and Smith \cite{HOFF} defined the concept of $internal\ path$ in a graph $G$ as a walk $v_0v_1 \cdots v_s$ $(s \geq 1)$ such that the vertices $v_0, v_1, \cdots ,v_s$ are all mutually distinct, $d(v_0) > 2$, $d(v_s) > 2$, and $d(v_i) = 2$ for $0 < i < s$.

\begin{clm}\label{cp}
If the graph $G$ contains at least two cycles, then it has at most one internal path.
\end{clm}

\noindent\textbf{Proof of Claim 4.}
If Claim 4 would not be true, then we could suppose that $G$ contains at least two internal paths, say $P$ and $P'$. Taking into account Claims 1-3, we conclude that such a graph $G$ is a chain cactus whose block is either a triangle or an edge and $G$ has no pendant path. That is to say, each of the end-block of $G$ is a triangle. Assume that the leftmost triangle-chain in $G$ is of length $i$ whose end-vertex is the end-vertex of the internal path $P$. Similarly, suppose that the rightmost triangle-chain in $G$ is of length $j$ whose end-vertex is an end-vertex of internal path $P'$. Without loss of generality, we assume that $j \geq i \geq 1$. By our assumption, it is routine to check that $i + j \leq t - 1$. Hence, one end-vertex of $P$ is the end-vertex of the leftmost triangle-chain of length $i$, the other one of $P$ is on a triangle, say $C^0$. This triangle is not in the rightmost (resp. leftmost) triangle-chain. Then we can partition $G$ as $G_1 \cup H \cup G_2$, where $G_1$ is composed of the leftmost $i$ successive triangle, $H$ is the internal path $P$ together with the triangle $C^0$, whereas $G_2$ is the rest subgraph of $G$. In other words, in the graph $G_2 \cong G[E(G) - E(G_1) \cup E(H)]$ is the connected component containing the internal path $P'$. It is straightforward to check that $| E(G_2) | \geq 3j + 1 > 3i = | E(G_1) |$. By Lemma~\ref{lem-11}, there exists another cactus graph $G'$ in $\mathcal{G}_{n,t}$, such that $W_e(G) < W_e(G')$, a contradiction to the choice of $G$.
\qed

\begin{clm}\label{cp}
The graph $G$ is isomorphic to the saw graph $Sw(i,j;n - 2i - 2j - 1)$ with $| i - j | \leq 1$.
\end{clm}

\noindent\textbf{Proof of Claim 5.}
In fact, $G \cong Sw(i,j;n - 2i - 2j - 1)$ follows directly from Claims 1-4. In order to complete the proof of Claim 5, it suffices to show that the saw graph $Sw(i,j;n - 2i - 2j - 1)$ satisfies $| i - j | \leq 1$. Otherwise, applying Lemma~\ref{lem-11} to $Sw(i,j;n - 2i - 2j - 1)$ yields another cactus graph $G'$ in $\mathcal{G}_{n,t}$ such that $W_e(G) < W_e(G')$, contradiction to the choice of $G$.
\qed

By Claims 1-5, we obtain that
$$G \cong Sw(\lfloor {t/2} \rfloor,\lceil {t/2} \rceil;n - 2t - 1).$$
By a direct calculation we obtain:
If $t = 2k$ with $k \geq 0$,
\[
W_e(Sw(\lfloor {t/2} \rfloor,\lceil {t/2} \rceil;n - 2t - 1)) = \frac{1}{6}n^3 - \frac{1}{2}n^2 + \frac{1}{3}n + kn^2 - 4k^2n + \frac{4}{3}k^3 + 8k^2 - \frac{10}{3}k.
\]
If $t = 2k + 1$ with $k \geq 0$,
\[
W_e(Sw(\lfloor {t/2} \rfloor,\lceil {t/2} \rceil;n - 2t - 1)) = \frac{1}{6}n^3 - \frac{13}{6}n + kn^2 - 4k^2n - 4kn - \frac{4}{3}k^3 + 10k^2 + \frac{35}{3}k + 5.
\]
\qed

\section{Concluding remarks}

In this paper, the lower bound on edge-Wiener index of the cacti with $n$ vertices and $t$ cycles is determined and the corresponding extremal graph is identified. Furthermore, the upper bound on edge-Wiener index of the cacti with given cycles is established and the corresponding extremal graph is given as well.
For further study, it would be interesting to determine the extremal graph that has the maximum edge Szeged index and revised edge Szeged index in these class of cacti.

\section*{Acknowledgments}
This work was supported by the National Natural Science Foundation of China
(Nos.  11731002, 11371052), the Fundamental Research Funds for the Central Universities (Nos. 2016JBM071, 2016JBZ012) and the $111$ Project of China (B16002).

\end{document}